\documentclass[12pt]{amsart}
\usepackage{amsmath, amssymb, latexsym, amsthm, mathrsfs, bm}
\usepackage[all]{xy}

\newcommand{\Pa}[9]{\bibitem{#1} {#2}, \emph{#3}, {#4} \textbf{#5} ({#6}), {#7}--{#8}.}

      \newenvironment{changemargin}[2]{\begin{list}{}{
         \setlength{\topsep}{0pt}\setlength{\leftmargin}{0pt}
         \setlength{\rightmargin}{0pt}
         \setlength{\listparindent}{\parindent}
         \setlength{\itemindent}{\parindent}
         \setlength{\parsep}{0pt plus 1pt}
         \addtolength{\leftmargin}{#1}\addtolength{\rightmargin}{#2}
         }\item }{\end{list}}
\newcommand{\znb}[2]{{\bm{[}#1; #2\bm{]}}}

\newcommand{\Setting}[7]{\xymatrix@R=4pt@C=7pt{#1\ar@{-}[r]&#2\ar@{-}[r]&#3\\&#4\ar@{-}[u]\\
#5\ar@{-}[uu]\ar@{-}[r] & #6\ar@{-}[u]\ar@{-}[r] & #7\ar@{-}[uu]}}

\newcommand{\Fr}{\mathit{F\!r}}

\newcommand{\compactN}{\cl{\mathbb{N}}}

\newcommand{\arx}[3]{\texttt{http://arxiv.org/#1/#3}}
\newcommand{\bq}{\begin{quote}}
\newcommand{\eq}{\end{quote}}
\newcommand{\cl}[1]{\overline{#1}}
\newcommand{\CH}{the Continuum Hypothesis}

\newcommand{\sr}[2]{{\txt{$#1$\\$#2$}}}
\newcommand{\N}{\mathbb{N}}
\newcommand{\NN}{{\N^{\N}}}

\newcommand{\NNup}{{\N^{\uparrow\N}}}

\newcommand{\roth}{{[\N]^{\aleph_0}}}

\newcommand{\Inc}{{\compactN^{\uparrow\N}}}

\newcommand{\NcompactN}{{\compactN^\N}}
\newcommand{\seq}[1]{\{#1\}_{n\in\N}}
\newcommand{\sseq}[1]{\{#1 : n\in\N\}}

\newcommand{\op}{\operatorname}

\newcommand{\cB}{\mathcal{B}}
\newcommand{\cC}{\mathcal{C}}

\newcommand{\CG}{C_\Gamma}

\newcommand{\cA}{\mathcal{A}}

\newcommand{\cM}{\mathcal{M}}

\newcommand{\cO}{\mathcal{O}}

\newcommand{\Q}{\mathbb{Q}}
\newcommand{\R}{\mathbb{R}}
\newcommand{\cU}{\mathcal{U}}

\newcommand{\cV}{\mathcal{V}}

\newcommand{\Impl}{\Rightarrow}
\long\def\forget#1\forgotten{}

\newcommand{\fb}{\mathfrak{b}}
\newcommand{\fc}{\mathfrak{c}}
\newcommand{\fd}{\mathfrak{d}}

\newcommand{\fp}{\mathfrak{p}}

\newcommand{\sbst}{\subseteq}

\newcommand{\sm}{\setminus}

\newcommand{\cov}{\op{cov}}

\newtheorem{thm}{Theorem}
\newcommand{\bthm}{\begin{thm}} \newcommand{\ethm}{\end{thm}}
\newtheorem{prop}[thm]{Proposition}
\newcommand{\bprp}{\begin{prop}} \newcommand{\eprp}{\end{prop}}
\newtheorem{fact}[thm]{Fact}
\newcommand{\bfct}{\begin{fact}} \newcommand{\efct}{\end{fact}}
\newtheorem{prob}[thm]{Problem}
\newcommand{\bprb}{\begin{prob}} \newcommand{\eprb}{\end{prob}}
\newtheorem{lem}[thm]{Lemma}
\newcommand{\blem}{\begin{lem}} \newcommand{\elem}{\end{lem}}
\newtheorem{cor}[thm]{Corollary}
\newcommand{\bcor}{\begin{cor}} \newcommand{\ecor}{\end{cor}}
\newtheorem{conj}[thm]{Conjecture}
\newcommand{\bcnj}{\begin{conj}} \newcommand{\ecnj}{\end{conj}}
\theoremstyle{definition}
\newtheorem{defn}[thm]{Definition}
\newcommand{\bdfn}{\begin{defn}} \newcommand{\edfn}{\end{defn}}
\theoremstyle{remark}
\newtheorem{rem}[thm]{Remark}
\newcommand{\brem}{\begin{rem}} \newcommand{\erem}{\end{rem}}
\newtheorem{exam}[thm]{Example}
\newcommand{\bexs}{\begin{exam}} \newcommand{\eexs}{\end{exam}}
\newcommand{\bpf}{\begin{proof}} \newcommand{\epf}{\end{proof}}

\newcommand{\be}{\begin{enumerate}}
\newcommand{\ee}{\end{enumerate}}
\newcommand{\bi}{\begin{itemize}}
\newcommand{\itm}{\item}
\newcommand{\ei}{\end{itemize}}

\forget 
 \forgotten

\newcommand{\sone}{\mathsf{S}_1}
\newcommand{\sfin}{\mathsf{S}_\mathrm{fin}}
\newcommand{\ufin}{\mathsf{U}_\mathrm{fin}}

\title{Continuous selections and $\sigma$-spaces}

\author{Du\v{s}an Repov\v{s}}
\address[Du\v{s}an Repov\v{s}]{Institute of  Mathematics, Physics and
Mechanics and Faculty of Education, University of Ljubljana,
P.O.B. 2964, Ljubljana, Slovenija 1001.}
\email{dusan.repovs@guest.arnes.si}

\author{Boaz Tsaban}
\address[Boaz Tsaban]{Department of Mathematics, Bar-Ilan University,
Ramat-Gan 52900, Israel;
and
Department of Mathematics,
Weizmann Institute of Science, Rehovot 76100, Israel.
}
\email{tsaban@math.biu.ac.il}

\author{Lyubomyr Zdomskyy}
\address[Lyubomyr Zdomskyy]{Department of Mathematics, Weizmann
Institute of Science, Rehovot 76100, Israel.}
\curraddr{Kurt G\"odel Research Center for Mathematical Logic, W\"ahringer Str.\ 25, A-1090 Vienna, Austria.}
\email{lzdomsky@gmail.com}

\thanks{The first and third authors were partially supported by the Slovenian Research Agency grants
P1-0292-0101-04 and BI-UA/04-06-007. The second author was partially supported by the Koshland Center for Basic Research.}



\begin{document}

\begin{abstract}
Assume that $X\sbst\R\sm\Q$, and each clopen-valued lower
semicontinuous multivalued map $\Phi:X\Rightarrow \Q$ has a
continuous selection $\phi:X \to \Q$. Our main result is that in
this case, $X$ is a $\sigma$-space. We also derive a partial
converse implication, and present a reformulation of the Scheepers
Conjecture in the language of continuous selections.
\end{abstract}

\maketitle

\section{Introduction}

All topological spaces considered in this note are assumed to have large inductive
dimension $0$, that is, disjoint closed sets can be separated by clopen sets.

By a \emph{multivalued map} $\Phi$ from a set $X$ into a set $Y$ we
understand a map from $X$ into  the power-set of
$Y$, denoted by $P(Y)$, and we write $\Phi:X\Rightarrow Y$.
Let $X,Y$ be topological spaces.
A multivalued map $\Phi:X\Rightarrow Y$ is \emph{lower semi-continuous (lsc)} if for each
open $V\sbst Y$, the set
$$\Phi^{-1}_\cap(V)=\{x\in X:\Phi(x)\cap V\neq\emptyset\}$$
is open in $X$.

A function $f:X\to Y$ is a \emph{selection} of a multivalued map $\Phi:X\Rightarrow Y$ if
$f(x)\in\Phi(x)$ for all $x\in X$. Let $\cC\sbst P(Y)$. A multivalued map $\Phi:X\Rightarrow Y$ is
\emph{$\cC$-valued} if $\Phi(x)\in\cC$ for all $x\in X$. Similarly,
we define \emph{clopen-valued}, \emph{closed-valued}, and \emph{open-valued}.
A general reference for selections of multivalued mappings is \cite{RS}.

\bthm[Michael \cite{Mi80}]\label{mich80}
Assume that $X$ is a countable space, $Y$ is a
first-countable space, and $\Phi:X\Rightarrow Y$ is lsc. Then $\Phi$ has
a continuous selection $\phi:X\to Y$.
\ethm

This result was extended in \cite[Theorem~3.1]{YJ}, where it was
proved that a space $X$ is countable if, and only if,
for each first-countable $Y$, each lsc multivalued map from $X$ to $Y$ has a continuous selection.
In fact, their proof gives the following.

\bthm[Yan-Jiang \cite{YJ}]\label{ch}
A separable space $X$ is countable if, and only if, for each first-countable space $Y$ and each
open-valued lsc map $\Phi:X\Rightarrow Y$, there is a continuous selection $\phi:X\to Y$.
\ethm

We extend Theorems \ref{mich80} and \ref{ch} by considering a qualitative restriction
on the space $X$ (instead of the quantitative restriction ``$X$ is countable'').
We also point out a connection to a conjecture of Scheepers.

\section{$\sigma$-spaces}

Define a topology on $P(\N)$ by identifying $P(\N)$ with the Cantor space $\{0,1\}^\N$.
The standard base of the topology of $P(\N)$ consists of the sets of the form
$$\znb{s}{t}=\{A\sbst\N :A\cap s=t\},$$
where $s$ and $t$ are finite subsets of $\N$.
Let $\Fr$ denote the \emph{Fr\'echet filter},
consisting of all cofinite subsets of $\N$,
and let $\roth$ be the family of all
infinite subsets of $\N$. $\Fr$ and $\roth$
are subspaces of $P(\N)$ and are homeomorphic to
$\Q$ and to $\R\sm\Q$, respectively (see \cite{Ke}).
Let
\begin{eqnarray*}
\cB & = & \{\znb{s}{\emptyset}:s \mbox{ is a finite subset of }\N\};\\
\cB_\Fr & = & \{B\cap \Fr:B\in\cB\}.
\end{eqnarray*}
Note that $\cB$ is the standard clopen base at the point $\emptyset\in P(\N)$.

A topological space $X$ is a \emph{$\sigma$-space} if each $F_\sigma$ subset of $X$ is a
$G_\delta$ subset of $X$ \cite{Mil}.

The main result of this note is the following.

\bthm\label{corref}
The following are equivalent:
\be
\itm $X$ is a $\sigma$-space;
\itm Each $\cB_\Fr$-valued lsc map $\Phi:X\Rightarrow\Fr$ has a continuous selection.
\ee
\ethm

The proof of Theorem \ref{corref} and subsequent results use the following notions.
A family $\cU=\{U_n:n\in\N\}$
of subsets of a set $X$ is a \emph{$\gamma$-cover} of $X$ if for
each $x\in X$, $x\in U_n$ for all but finitely many $n$.
A bijectively enumerated family $\cU=\{U_n:n\in\N\}$ of subsets of a set $X$
induces a \emph{Marczewski map} $\cU:X\to P(\N)$ defined by
$$\cU(x)=\{n\in\N : x\in U_n\}$$
for each $x\in X$ \cite{M-S}.

\begin{rem}
Marczewski maps can be naturally associated to any sequence of
sets, not necessarily bijectively enumerated. Our restriction to bijective
enumerations allows working with the classical notion of $\gamma$-cover.
An alternative approach would be to use \emph{indexed} $\gamma$-covers, that
is, sequences of sets $(U_n : n\in\N)$ such that each $x\in X$ belongs
to $U_n$ for all but finitely many $n$. All results of the present paper
hold in this setting, too.
\end{rem}

For a function $f:X\to Y$, $f[X]$ denotes $\{f(x) : x\in X\}$, the image of $f$.

\blem \label{gamma}
Let $\cU=\{U_n:n\in\N\}$ be a bijectively enumerated family of subsets of a
topological space $X$. Then
\be
\item $\cU$ is a clopen $\gamma$-cover of $X$ if, and only if,
$\cU[X]\sbst\Fr$ and $\cU:X\to P(\N)$ is continuous;
\item $\cU$ is an open $\gamma$-cover
of $X$ if, and only if, $\cU[X]\sbst\Fr$ and the multivalued map
$\Phi:X\Rightarrow\Fr$ defined by $\Phi(x)=P(\cU(x))\cap \Fr$ is lsc.
\ee
\elem
\bpf
The first assertion follows immediately from the corresponding
definitions. To prove the second assertion, let us assume that
$\cU=\{U_n:n\in\N\}$ is an open $\gamma$-cover of $X$. Fix some
finite subsets $s,t$ of $\N$ and $x\in X$ such that
$\znb{s}{t}\cap\Phi(x)\neq\emptyset.$ There exists $A\in \Fr$ such
that $A\sbst\cU(x)$ and $A\cap s=t$.

Let $V=\bigcap_{n\in\cU(x)\cap s} U_n$. The set $V$ is
open in $X$, being an intersection of finitely many open sets, and it
contains $x$ by definition of $\cU$.
Thus it suffices to show that $\znb{s}{t}\cap\Phi(y)\neq\emptyset$
for all $y\in V$. A direct verification indeed
shows that $(A\cap s)\cup (\cU(y)\setminus s)$ belongs to $\znb{s}{t}$ as well as to
$\Phi(y)$.

To prove the converse implication, it suffices to note that
$U_n=\Phi_\cap^{-1}\: \znb{\{n\}}{\{n\}}$, and use the lower
semi-continuity of  $\Phi$.
\epf

The following is a key result of Sakai.
A cover $\{U_n:n\in\N\}$ of $X$ is \emph{$\gamma$-shrinkable}
\cite{Sa} if there is a clopen $\gamma$-cover $\{C_n:n\in\N\}$ of $X$
such that $C_n\sbst U_n$ for all $n$.
Note that $\cU$ is a $\gamma$-cover of $X$ if, and only if, $\cU[X]\sbst\Fr$.

\bthm[Sakai \cite{Sa}]\label{sakai}
$X$ is a $\sigma$-space if, and only if, each open $\gamma$-cover of $X$ is $\gamma$-shrinkable.
\ethm

\bpf[Proof of Theorem \ref{corref}]
$(2\Impl 1)$ Assume that each $\cB_\Fr$-valued lsc $\Phi:X\Rightarrow \Fr$
has a continuous selection. We will show that $X$ is a $\sigma$-set
by using Sakai's characterization (Theorem \ref{sakai}).

Let $\cU$ be an open $\gamma$-cover of $X$.
Define $\Phi(x)=P(\cU(x))\cap\Fr$. $\Phi$ is $\cB_\Fr$-valued, and by Lemma \ref{gamma},
$\Phi$ is lsc. By our assumption, $\Phi$ has a continuous selection.
The following lemma implies that $\cU$ is $\gamma$-shrinkable.

\blem \label{gamref}
Let $\cU=\{U_n:n\in\N\}$ be a bijectively enumerated open $\gamma$-cover of a space $X$.
The following are equivalent:
\be
\itm $\cU$ is $\gamma$-shrinkable;
\itm The multivalued map $\Phi(x)=P(\cU(x))\cap \Fr$ has a
continuous selection.
\ee
\elem
\bpf
$(1\Impl 2)$ If $\mathcal V=\sseq{V_n}$ is a witness for $(1)$, then the
map $x\mapsto\cV(x)$ is a continuous selection of $\Phi$.

$(2\Impl 1)$ If $\phi:X\to\Fr$ is a
continuous selection of $\Phi$, then $\{V_n:=\{x\in X:\phi(x)\ni
n\}:\ n\in\N\}$ is a clopen $\gamma$-cover of $X$ with the
property $V_n\sbst U_n$, for all $n\in\N$. Indeed, if $x\in V_n$,
then $n\in\phi(x)\in P(\cU(x))\cap\Fr$, and hence
$n\in\cU(x)$, which is equivalent to $x\in U_n$.
\epf

$(1\Impl 2)$ Assume that $X$ is a $\sigma$-space and
$\Phi:X\Rightarrow\Fr$ is lsc and $\cB_\Fr$-valued.
The following is easy to verify.

\blem \label{obs11}
For each $\cB_\Fr$-valued $\Phi:X\Rightarrow\Fr$, there exists a map $\phi:X\to\Fr$ such that
$\Phi(x)=P(\phi(x))\cap\Fr$ for all $x\in X$.

Conversely, for each map $\phi:X\to\Fr$, the multivalued map
$\Phi:X\Rightarrow\Fr$ defined by $\Phi(x)=P(\phi(x))\cap\Fr$
is $\cB_\Fr$-valued.\hfill\qed
\elem

Let $\phi$ be as in Lemma~\ref{obs11}. For each $n$, let
$U_n=\{x\in X:n\in\phi(x)\}=\{x\in X:\Phi(x)\cap \znb{\{n\}}{\{n\}}\neq\emptyset\}$.
Each $U_n$ is open, and $\cU=\sseq{U_n}$ is a $\gamma$-cover of $X$.
By Sakai's Theorem \ref{sakai}, $\cU$ is $\gamma$-shrinkable.

Note that the Marczewski map induced by the family $\cU$ is exactly the map $\phi$.
Thus, by Lemma \ref{gamref}, $\Phi(x)=P(\phi(x))\cap\Fr=P(\cU(x))\cap\Fr$ has a continuous selection.
\epf

\bcor
If each clopen-valued lsc map $\Phi:X\Rightarrow\Fr$ has a continuous selection,
then $X$ is a $\sigma$-space.\hfill\qed
\ecor

\bprb
Assume that $X\sbst\R$ is a $\sigma$-space.
Does each clopen-valued lsc map $\Phi:X\Rightarrow\Q$ have a continuous selection?
\eprb

It is consistent (relative to ZFC) that each metrizable separable
$\sigma$-space $X$ is countable \cite{Mil}.
Thus, by Theorems~\ref{mich80} and \ref{corref},
we have the following extension of Theorem~\ref{ch}.

\bcor \label{cor1}
It is consistent that the following statements
are equivalent, for metrizable separable spaces $X$:
\be
\item Every clopen-valued lsc map $\Phi:X\Rightarrow\Q$ has a continuous selection;
\item $X$ is countable.\hfill\qed
\ee
\ecor

\bprb
Is Corollary~\ref{cor1} provable in ZFC?
\eprb

$\fb$ is the minimal cardinality
of a subset of $\NN$ which is unbounded with respect to $\leq^*$
($f\le^* g$ means: $f(n)\le g(n)$ for all but finitely many $n$.)
$\fb$ is uncountable, and consistently, $\aleph_1<\fb$ \cite{BlassHBK}.
If $|X|<\fb$, then $X$ is a $\sigma$-set \cite{FM88, CBC}. By Theorem
\ref{corref}, we have the following quantitative result.

\bcor
Assume that $|X|<\fb$.
Then for each $\cB_\Fr$-valued lsc map $\Phi:X\Rightarrow\Fr$, $\Phi$ has a continuous selection.\hfill\qed
\ecor

\section{$\fb$-scales}

Let $\NNup$ be the set of all (strictly) increasing elements of
$\NN$.
$B=\{b_\alpha:\alpha<\fb\}\sbst \NNup$ is a \emph{$\fb$-scale}
if $b_\alpha\leq^* b_\beta$ for all $\alpha<\beta$, and $B$ is unbounded with respect
to $\leq^*$.
$\compactN=\N\cup\{\infty\}$ is a convergent sequence with the limit point $\infty$,
which is assumed to be larger than all elements of $\N$.
$\Inc$ is the set of all nondecreasing elements of $\NcompactN$, and
$Q=\{x\in\Inc:(\exists m)(\forall n\ge m)\ x(n)=\infty\}$ is the set of all
``eventually infinite'' elements of $\Inc$.

Sets of the form $B\cup Q$ where $B$ is a $\fb$-scale
were extensively studied in the literature (see \cite{BT,Mil,TZ} and references therein).
$B\cup Q$ is concentrated on $Q$ and is therefore not a $\sigma$-space.
Consequently, it does not have the properties stated in Theorem~\ref{corref}.
In fact, we have following.

\bthm\label{1.7}
Let $X=B\cup Q$, where $B\sbst\NN$ is a $\fb$-scale. Then there
exists a clopen-valued lsc map
$\Phi:X\Rightarrow\Q$ with the following properties:
\be
\item[$(i)$]  $\Phi(x)=\Q$, for all $x\in B$; and
\item[$(ii)$] For each $Y\sbst X$ such that $Q\sbst Y$,
and each continuous $\phi:Y\to\Q$ such that $\phi(y)\in\Phi(y)$ for all $y\in Y$,
$|Y|<|X|$.
\ee
\ethm
\bpf
Write $Q=\{q_n:n\in\N\}$, and consider the $\gamma$-cover
$\cU=\{U_n:n\in\N\}$ of $X$, where
$U_n=X\setminus\{q_n\}$, $n\in\N$.

\blem\label{cl1.7}
For each $B'\sbst B$ with $|B'|=\fb$,
and each choice of clopen sets $V_n\sbst U_n$, $n\in\N$,
there is $x\in B'$ such that $\{n : x\notin V_n\}$ is infinite.
\elem
\bpf
Assuming the converse, we could find a clopen
$\gamma$-cover $\{V'_n:n\in\N\}$ of $B'\cup Q$ such that
$V'_n\sbst U_n$. Let $V_n$ be a closed subspace of $\Inc$ such that
$V_n\cap X=V'_n$. Then $W_n=\Inc\setminus V_n$ is an open
neighborhood of $q_n$ in $\Inc$. Set $G_n=\bigcup_{k\geq n}W_k$ and
$G=\bigcap_{n\in\N}G_n$. For each $n\in\N$ the set $\NNup\cap
(\Inc\setminus G_n)$ is  a cofinite subset of the compact space
$\Inc\setminus G_n$, and hence it is $\sigma$-compact.

Therefore
$\NNup\cap (\Inc\setminus G)=\bigcup_{n\in\N}\NNup\cap (\Inc\setminus
G_n)$ is a $\sigma$-compact subset of $\NNup$ as well. Since $B'$
is unbounded, there exists  $x\in B'\cap G$, and hence $x$ belongs
to $W_n$ for infinitely many $n\in\N$, which implies that
$\{n\in\N:x\not\in V_n\}=\{n\in\N:x\not\in V'_n\}$ is infinite, a
contradiction.
\epf

Recall that $\Fr$ is homeomorphic to $\Q$. Thus, it suffices to construct a lsc
 $\Psi:X\Rightarrow\Fr$ with the  properties $(i)$
and $(ii)$. Set $\Psi(x)=P(\cU(x))\cap\Fr$.

By Lemma~\ref{gamma}, the multivalued map $\Psi$ is lsc  and
there are no partial selections $f:Y\to \Fr$ defined on subsets
$Y\sbst X$ such that $|Y|=|X|=\fb$ and $Y\supset Q$. Indeed, it
suffices to use Lemmata \ref{gamref} and \ref{cl1.7},
asserting that  there is no clopen refinement $\{V_n:n\in\N\}$ of
$\{U_n:n\in\N\}$ which is a $\gamma$-cover of such a subspace $Y$
of $X$.
\epf

Theorem \ref{1.7} can be compared with Theorem~1.7 and Example~9.4 of \cite{Mi80}.

\medskip

The undefined terminology in the following discussion is standard and can be found in,
e.g., \cite{Sc}. Lemma \ref{cl1.7} motivates the
introduction of the following covering property of a space $X$:
\begin{description}
\itm[$(\theta)$]
There exists an open $\gamma$-cover $\cU=\{U_n:n \in
\N\}$ of  $X$ and a countable  $D\sbst X$ such that for any family
$\cV=\{V_n:n \in \N\}$ of clopen subsets of $X$ with $V_n\sbst U_n$
for all $n$,  if $\cV$ is a $\gamma$-cover of some $Y\sbst X$ such
that $D\sbst Y$, then $|Y|<|X|$.
\end{description}
Theorem \ref{1.7} implies the following.

\begin{cor}\label{happ}
Assume that $X=B\cup Q$ where $B\sbst\NN$ is a $\fb$-scale. Then
$X$ satisfies $(\theta)$.\hfill\qed
\end{cor}

\begin{figure}[!htp]
{
\renewcommand{\sr}[2]{{#1}}
\begin{changemargin}{-4cm}{-3cm}
\begin{center}
$\xymatrix@R=10pt{
&
&
& \sr{\ufin(\cO,\Gamma)}{\fb}\ar[r]
& \sr{\ufin(\cO,\Omega)}{\fd}\ar[rr]
& & \sr{\sfin(\cO,\cO)}{\fd}
\\
&
&
& \sr{\sfin(\Gamma,\Omega)}{\fd}\ar[ur]
\\
& \sr{\sone(\Gamma,\Gamma)}{\fb}\ar[r]\ar[uurr]
& \sr{\sone(\Gamma,\Omega)}{\fd}\ar[rr]\ar[ur]
& & \sr{\sone(\Gamma,\cO)}{\fd}\ar[uurr]
\\
&
&
& \sr{\sfin(\Omega,\Omega)}{\fd}\ar'[u][uu]
\\
& \sr{\sone(\Omega,\Gamma)}{\fp}\ar[r]\ar[uu]
& \sr{\sone(\Omega,\Omega)}{\cov(\cM)}\ar[uu]\ar[rr]\ar[ur]
& & \sr{\sone(\cO,\cO)}{\cov(\cM)}\ar[uu]
}$
\end{center}
\end{changemargin}
}
\caption{The Scheepers Diagram}\label{schdiag}
\end{figure}
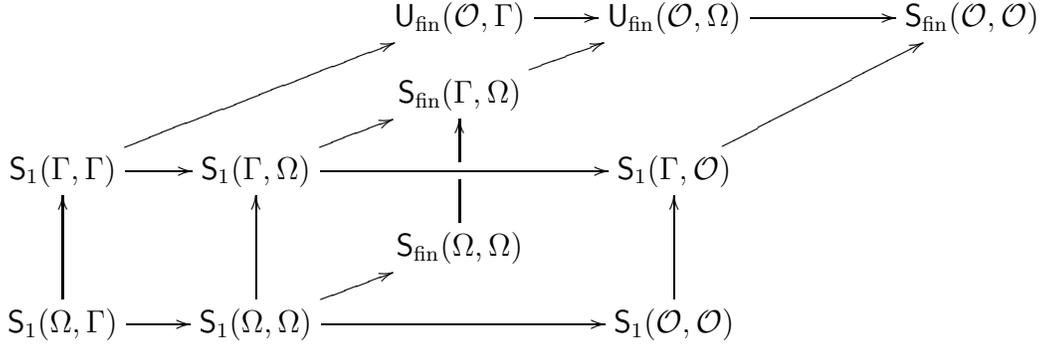

The property $(\theta)$ seems to stand apart  from
the classical selection principles considered in \cite{Sc, coc2}.
Figure \ref{schdiag} (reproducing Figure 3 on page 245 of \cite{coc2})
summarizes the relations among these properties.

Every countable space satisfies the strongest property in that
figure, namely $\sone(\Omega,\Gamma)$ \cite{GM},
and it is clear that countable spaces do not satisfy $(\theta)$.
Moreover, by Sakai's Theorem \ref{sakai}, no $\sigma$-space
satisfies $(\theta)$.

Assuming \CH{} there is a $\fb$-scale $B$ such that
$B\cup Q$ is not a $\sigma$-space, but satisfies $\sone(\Omega,\Gamma)$ \cite{GM}
as well as $(\theta)$ (Corollary \ref{happ}).

Consider the topological sum $X=\R\oplus(\R\sm\Q)$.
The open sets $U_n=(-n,n)\oplus(\R\sm\Q)$, $n\in\N$, form a $\gamma$-cover of $X$
and show that $X$ satisfies $(\theta)$ for a trivial reason, and
does not satisfy the weakest property in the Scheepers Diagram, namely $\sfin(\cO,\cO)$,
because it contains $(\R\sm\Q)$ as a closed subspace.
A less trivial (zero-dimensional) example is given in the following consistency result.

\bthm
Assume that $\fb=\fd=\mathrm{cf}(\fc)<\fc$.
There is a set $X\sbst\R\sm\Q$ satisfying $(\theta)$
but not $\sfin(\cO,\cO)$.
\ethm
\bpf
Let $B=\{b_\alpha:\alpha<\fb\}$ be a $\fb$-scale and
$\fc=\bigcup_{\alpha<\fb}\lambda_\alpha$ with
$\lambda_\alpha<\fc$. Fix $D_\alpha\sbst\NNup$ such that
$|D_\alpha|=\lambda_\alpha$ and for each $f\in D_\alpha$, $|f(n)-b_\alpha(n)|<2$ for all $n$.

Let $Y\sbst\NN$ be a dominating family.
The direct sum of $X=Q\cup\bigcup_{\alpha<\fb}D_\alpha$ and $Y$
satisfies $(\theta)$ by the methods of Theorem~\ref{1.7}.
But $Y$ is a closed subset of this space and does not satisfy $\sfin(\cO,\cO)$ \cite{Sc}.
\epf

\section{Connections with the Scheepers Conjecture}

Let $\cA$ and $\cB$ be any two families. Motivated by works of
Rothberger, Scheepers introduced the following prototype of
properties \cite{Sc}:
\begin{description}
\item[$\sone(\cA,\cB)$]
For each sequence $\seq{\cU_n}$ of members of $\cA$, there
exist members $U_n\in \cU_n$, $n\in\N$, such that
$\{U_n:n\in\N\}\in\cB$.
\end{description}
Let $\Gamma$ and $C_\Gamma$ be the collections of all open and clopen $\gamma$-covers of a
set $X\sbst\R$, respectively. Scheepers  \cite{Sc1}
has conjectured that  the property $\sone(\Gamma,\Gamma)$ is equivalent to a certain local property
in the space of continuous real-valued functions on $X$.
Sakai \cite{Sa} and independently Bukovsk\'y-Hale\v{s} \cite{BH07}
proved that Scheepers' Conjecture holds if, and only if,
$\sone(\Gamma,\Gamma)=\sone(\CG,\CG)$ for sets of reals.

Lemma~\ref{gamma} establishes a bijective
correspondence between open $\gamma$-covers of a space $X$ and
maps $\phi:X\to \Fr$ for which the multivalued map
$\Phi(x)=P(\phi(x))\cap \Fr$ is lsc. This is used
in the proof of the following characterizations, which
give an alternative justification for the Scheepers Conjecture.

\bthm\label{hi}
$X$ satisfies $\sone(\CG,\CG)$ if, and only if,
for each continuous $\phi:X\to\Fr^\N$ there is
$f\in\NN$ such that $f(k)\in\phi(x)(k)$ for each
$x\in X$ and all but finitely many $k$.
\ethm

Since the proof of Theorem \ref{hi} is easier than that of
the following theorem, we omit it.

\bthm
$X$ satisfies $\sone(\Gamma,\Gamma)$ if, and only if,
for each $\phi:X\to\Fr^\N$ such that the multivalued map
$\Phi:x\mapsto\Pi_{k\in\N}(P(\phi(x)(k))\cap\Fr)$ is
lsc, there is $f\in\NN$ such that
$f(k)\in\phi(x)(k)$ for each $x\in X$ and all but finitely many
$k$.
\ethm
\bpf
Assume that $X$ satisfies $\sone(\Gamma,\Gamma)$.
Fix a map $\phi:X\to\Fr^\N$ as in the second assertion. The multivalued map
$\Phi_k:X\Rightarrow\Fr$ assigning to each point $x\in X$ the subset
$\Phi_k(x)=P(\phi(x)(k))\cap\Fr$ of $\Fr$,
is lsc for all $k$.

The family $\{U_{k,n}:n\in\N\}$, where
$U_{k,n}=\{x\in X:\Phi_k(x)\cap
\znb{\{n\}}{\{n\}}\neq\emptyset\}=\{x\in X:n\in\phi(x)(k)\}$, is an
open $\gamma$-cover of $X$. Since $X$ satisfies $\sone(\Gamma,\Gamma)$,
there exists $f\in\NN$ such that
$\{U_{k,f(k)}:k\in\N\}$ is a $\gamma$-cover of $X$. This implies
that $f(k)\in\phi(x)(k)$ for all $x\in\N$ and all but finitely many
$k$.

The proof of the converse implication is similar, using Lemma~\ref{gamma}.
\epf

\subsection*{Acknowledgements} We thank Taras Banakh for stimulating conversations.

\end{document}